\documentclass[preprint,11pt,natbib,fleqn,authoryear]{elsarticle}
        
\usepackage{bm}

\usepackage{amsmath}
\usepackage{amssymb}
\usepackage{theorem}
\usepackage{xspace}
\usepackage{geometry}
\usepackage{array}
\usepackage{subfigure}
\usepackage{epsfig}
\usepackage{hhline}
\usepackage{bbm}
 
\usepackage[utf8]{inputenc}    
\usepackage[francais]{babel}
\DeclareMathAccent{\widehat}{\mathord}{largesymbols}{"62}
\DeclareMathAccent{\widetilde}{\mathord}{largesymbols}{"65}
\def\pth#1{\left(#1\right)}
\def\acc#1{\left\{#1\right\}}

\def\eeX{\mathbb{X}}

\def\ebo{\textrm{\mathversion{bold}$\mathbf{\beta^0}$\mathversion{normal}}}

\def\eb{\textrm{\mathversion{bold}$\mathbf{\beta}$\mathversion{normal}}}  
\def\ed{\textrm{\mathversion{bold}$\mathbf{\delta}$\mathversion{normal}}}

\def\eU{\textrm{\mathversion{bold}$\mathbf{\Upsilon}$\mathversion{normal}}} 

\def\eE{\mathbb{E}}

\def\e1{1\!\!1}

\def\eeX{\mathbb{X}}

\theoremstyle{plain}

\newtheorem{theorem}{Theorem}[section]

\newtheorem{remark}{Remark}[section]

\newcommand{\beqn}{\begin{eqnarray*}}
\newcommand{\eeqn}{\end{eqnarray*}}
  
\def\ee1{\textrm{\mathversion{bold}$\mathbf{\varepsilon}$\mathversion{normal}}}

\def\eu{\mathbf{{u}}}

\newcommand{\N}{\mathbb{N}}
\newcommand{\R}{\mathbb{R}}
\newcommand{\PP}{\mathbb{P}}

\def\eX{\mathbf{X}}

\newcommand{\Var}{\mathbb{V}\mbox{ar}\,}

\def\argmin{\mathop{\mathrm{arg\,min}}} 
\def\hh{ \hspace*{0.5cm}}
\begin{document}
\begin{frontmatter}
\title {{\bf Variable selection  in high-dimensional linear model with possibly asymmetric or heavy-tailed   errors}}
 
\author[label3]{  Gabriela CIUPERCA}



 \address{Institut Camille Jordan, Université Lyon 1, France\fnref{label3}}

\fntext[label3]{Université de Lyon, Université Lyon 1, CNRS, UMR 5208, Institut Camille Jordan, Bat.  Braconnier, 43, blvd du 11 novembre 1918, F - 69622 Villeurbanne Cedex, France\\
\textit{Email address}: Gabriela.Ciuperca@univ-lyon1.fr}
 
 \begin{abstract}
We consider   the problem of automatic variable selection in a linear model with asymmetric or heavy-tailed errors  when the number of explanatory variables diverges with the sample size. For this high-dimensional model, the penalized least square method is not appropriate and the quantile framework makes the inference more difficult because to the non differentiability of the loss function. We propose and study an estimation method by penalizing the expectile process with an adaptive LASSO penalty. Two cases are considered: the number of model parameters is   smaller and afterwards larger  than the sample size, the two cases being distinct by the adaptive penalties considered. For each case we give the rate convergence and  establish the oracle properties of the adaptive LASSO expectile estimator. The proposed estimators are evaluated through Monte Carlo simulations and compared with the adaptive LASSO quantile estimator. We applied also our estimation method to real data in genetics when the number of parameters is greater than the sample size. 
 \end{abstract}
\end{frontmatter}

\section{Introduction}
The focus of the present paper is to better detect significant variables in a linear model, with the possibility that the number of explanatory variables is greater than the number of observations and when the error distribution is  asymmetrical or heavy-tailed. For this type of law, the using of the least squares (LS) estimation method is not appropriate for the  estimator accuracy. One possibility would be to use the quantile method, but it has the disadvantage that the loss function is not derivable, which complicates the theoretical study  but also computational methods.  A very interesting possibility is to consider the expectile method, introduced by \cite{Newey-Powell.87}, under assumption that the first moments of  $\varepsilon$ exist. This method has the advantage that the loss function is differentiable, which simplifies the theoretical study  and facilitates the numerical calculation. \\
In application fields (genetics, chemistry, biology, industry, finance), with the development in recent years of storage and/or measurement tools, we are confronted to study the influence of a very large number of variables on a studied process. That is why, let us consider in the present work the following linear model:  
\begin{equation}
\label{eq1}
Y_i=\eX_i^t \eb+\varepsilon_i, \qquad i=1, \cdots , n,
\end{equation}
with the vector parameter $\eb \in \R^p$ and $\ebo$ its true value (unknown). The size $p$ of $\ebo$ can depend on $n$ but the components $\eb^0_j$ don't depend on $n$, for any  $j=1, \cdots , p$. The vector $\eX_i=(X_{i1}, \cdots, X_{ip})$ contains the values of the  $p$  explanatory deterministic variables and $Y_i$ the values of response variable for observation $i$. The values $(Y_i, \eX_i)$ are known for any $i=1, \cdots , n $. Throughout the paper, all vectors are column.  If $p$ is very large, in order to find the explanatory variables that significantly influence the response variable  $Y$, an automatic selection should be made without performing hypothesis tests. When $p$ is very large, the use of hypothesis tests is not appropriate because they provide results with great variability and the significant explanatory variable final choice is not optimal (see \cite{Breiman-96}). \\
For model (\ref{eq1}), let then  the index set of the non-null true parameters, 
\[
{\cal A} \equiv \big\{   j \in \{ 1, \cdots, p\}; \, \beta^0_j \neq 0\big\}.
\]
Since $\ebo$ is unknown then, the set ${\cal A}$ is also unknown. We assume, without reducing generality, that   ${\cal A}=\{1, \cdots, p_0 \}$ and its complementary set is  ${\cal A}^c=\{p_0+1, \cdots , p  \}$, with $p_0 \leq p$. Hence the first $p_0$ explanatory variables have a significant influence on the response variable and the last $p-p_0$ variables are irrelevant.   Thus, the true parameter vector can be written as  $\ebo=\big(\eb^0_{\cal A}, \eb^0_{{\cal A}^c}\big)=\big(\eb^0_{\cal A},\textbf{0}_{p-p_0}\big)$, with $\textbf{0}_{p-p_0}$ a $(p-p_0)$-vector with all components zero. The number $p_0$ of the nonzero coefficients can depends on $n$. \\
For a vector $\eb$ we use the notational convention   $\eb_{\cal A}$ for its subvector containing the corresponding components of ${\cal A}$. For $i=1, \cdots, n$, we denote by $ \eX_{i,{\cal A}} $ the $p_0$-vector with the components $X_{ij}$, $j=1, \cdots, p_0$. We also use the notation $|{\cal A}|$ or $Card({\cal A})$ for the cardinality of ${\cal A}$. \\
In order to find the elements of ${\cal A}$, one of the most used techniques is the adaptive LASSO method, introduced by \cite{Zou.06} for  $p$ fixed by penalizing the squares sum  with an weighted $L_1$ penalty. This type of  parameter estimator is interesting  if it satisfies the oracle properties, i.e. the two following properties occur: 
\begin{itemize}
 \item \textit{sparsity of estimation}: the non-zero parameters are estimated as non-zero and the null parameters are shrunk directly as 0, with a probability converging to 1 when $n \rightarrow \infty$;
 \item \textit{asymptotic normality} of non-zero parameter estimators. 
 \end{itemize} 
 In order to distinguish between different types of adaptive LASSO estimators, we will use the term "adaptive LASSO LS-estimator" to refer to the minimizer of the LS sum penalized with adaptive LASSO. \\
Give some papers from very rich literature that consider the adaptive LASSO LS-estimator when $p$ depends on $n$, with the possibility that  $p>n$ are: \cite{Huang-Ma-Zhang.08}, \cite{Wang-Kulasekera.12},  \cite{Yang-Wu.16}. If the moments of the errors do not exist or the distribution of $\varepsilon$ presents outliers, then the  LS framework is not appropriate. One possibility is to consider the quantile model with the adaptive LASSO penalty. The recent literature is also very rich: \cite{FFB14}, \cite{Kaul-Koul.15}, \cite{Tang-Song-Wang-Zhu.13}, \cite{Ciuperca-18}, \cite{Zheng-Gallagher-Kulasekera-13},  \cite{Zheng-Peng-He.15}, to give just a few examples. As stated before, the loss function for quantile method not differentiable complicates the theoretical study and its computational implementation, which is a very important aspect in high-dimensionality.
Let us mention another work of \cite{Fan-Li-Wang.17} which is also devoted to the high dimensional regression in absence of symmetry of the errors and which proposes a penalized Huber loss. Only that the Huber loss function is not differentiable. 
 From where the idea of considering the expectile loss function for an high-dimensional model. In order to introduce the expectile method, for a fixed $\tau \in (0,1)$, let us consider the function $\rho_\tau(.)$ of the form 
\[\rho_\tau(x)=|\tau - \e1_{x <0}|x^2, \qquad \textrm{with} \quad x \in \R.
\]
\hh For the error and the design of model (\ref{eq1}) we make the following basic assumptions. \\
 The errors $\varepsilon_i$  satisfy the following  assumption:\\
 \textbf{(A1)} $(\varepsilon_i)_{1 \leqslant i \leqslant n}$ are i.i.d. such that $\eE[\varepsilon^4_i]< \infty$ and $
\eE[\varepsilon\big(\tau \e1_{\varepsilon >0}+ (1-\tau)\e1_{\varepsilon<0} \big)]=0$, that is its  $\tau$-th expectile is zero: $\eE[\rho'_\tau(\varepsilon)]=0$.\\
While, the design $(\eX_i)_{1 \leqslant i \leqslant n}$ satisfies the following assumption:\\
\textbf{(A2)} there exists two positive constants $m_0, M_0$ such that, $0 < m_0 \leq \mu_{min}\big(n^{-1} \sum^n_{i=1} \eX_i \eX_i^t\big)$ $\leq \mu_{max}\big(n^{-1} \sum^n_{i=1} \eX_i \eX_i^t\big) \leq M_0 < \infty$.\\

For a positive definite matrix, we denote by $\mu_{min}(.)$ and $\mu_{max}(.)$ its largest and smallest eigenvalues, respectively. Let us consider $\varepsilon$ the generic variable for the sequence $(\varepsilon_i)_{1 \leqslant i \leqslant n}$.  
Assumption (A1) is commonly required for the expectile models, see  \cite{Zhao-Chen-Zhang.18}, \cite{Gu-Zou.16}, \cite{Liao-Park-Choi.18}, \cite{Newey-Powell.87}, while assumption (A2) is standard in linear model for the parameter identifiability (considered also by \cite{Gao-Huang.10}, \cite{Wang-Wang.14}, \cite{Fan-Li-Wang.17}, \cite{Zou-Zhang-09}). 
Other assumptions will be stated about design in the following two sections, depending of the size  $p$ which varies in turn with $n$.\\
Quite in general, it is wise to use the expectile method when the moments of $\varepsilon$ exist but its distribution is  asymmetric or heavy-tailed. For $\tau=0.5$, we get the classical method of least squares. \\
For model (\ref{eq1}), consider the expectile process
\[
Q_n(\eb)\equiv \sum^n_{i=1} \rho_\tau (Y_i - \eX_i^t \eb)
\]
and the one with LASSO adaptive penalty:
\begin{equation}
\label{Rnb}
R_n(\eb)\equiv \sum^n_{i=1} \rho_\tau (Y_i - \eX_i^t \eb)+n \lambda_n \sum^p_{j=1} \widehat \omega_{n,j} |\beta_j|.
\end{equation}
The adaptive weights $ \widehat \omega_{n,j}$  will be defined later depending on whether  $p$  is smaller or larger than   $n$. The tuning parameter $\lambda_n$ is a positive deterministic sequence which together  $\widehat \omega_{n,j}$ controls the overall model complexity. Hence, we should choose $\lambda_n$ and $\widehat \omega_{n,j}$ such that $n \lambda_n \widehat \omega_{n,j}  \overset{\PP} {\underset{n \rightarrow \infty}{\longrightarrow}} 0$ for non-null parameters and   $n \lambda_n \widehat \omega_{n,j}  \overset{\PP} {\underset{n \rightarrow \infty}{\longrightarrow}} \infty $ for null coefficients. 
In order to automatically detect the null and non-zero components of $\eb$, we proceed in a similar way as for the adaptive LASSO LS estimation introduced by \cite{Zou.06}, and we consider the adaptive LASSO expectile estimator of $\eb$:
\begin{equation}
\label{hbetan}
 \widehat \eb_n \equiv \argmin_{\eb \in \R} R_n(\eb).
\end{equation}
The components of $\widehat \eb_n $ are $\widehat \eb_n =\big( \widehat{\beta}_{n,1}, \cdots , \widehat{\beta}_{n,p} \big)$. Similarly to  ${\cal A}$, let's define the index set:
\[
\widehat{{\cal A}}_n \equiv \{ j \in \{1, \cdots , p \}; \; \widehat{\beta}_{n,j} \neq 0 \},
\]
with the non-zero components of the adaptive LASSO expectile estimator.\\
The estimator $\widehat \eb_n$ will satisfy the oracle properties if:
\begin{itemize}
\item \textit{sparsity}: $\lim_{n \rightarrow \infty} \PP\big[{\cal A} = \widehat{{\cal A}}_n \big]=1$.
\item \textit{asymptotic normality}: for any vector $\eu \in \R^{p_0}$ with bounded norm, we have that:\\  $\sqrt{n} (\eu^t \eU^{-1}_{n,{\cal A}} \eu)^{-1/2} \eu^t ( \widehat{\eb}_n - \eb^0)_{\cal{A}} $ converges in distribution to a zero-mean   Gaussian law, with the $p_0$-squared matrix: $\eU_{n,{\cal A}} \equiv n^{-1} \sum^n_{i=1} \eX_{i,{\cal A}}  \eX_{i,{\cal A}}^t$.  
\end{itemize}
The asymptotic distribution of the parameter estimators can be used to construct asymptotic hypothesis tests and confidence intervals for the non-null parameters. 
For $p$ fixed, the properties of the estimator $\widehat \eb_n$ have been studied by \cite{Liao-Park-Choi.18} where it is shown that the convergence rate of $\widehat{\eb}_n$   towards $\ebo$ is of order  $n^{-1/2}$ and  that  $\widehat \eb_n$ satisfies the oracle properties. The case $p$ fixed was also studied by \cite{Zhao-Zhang.18} which consider a  penalized linear expectile regression with SCAD penalty function and obtain a  $n^{-1/2}$-consistent  estimator with oracle properties.  In the present paper we consider   $p$ depends on $n$, more precisely, of the form $p=O(n^c)$, with the constant $c >0$. The size $p_0$ and    the set  ${\cal A}$ can depend on $n$.  
The case when $p$ depends on $n$ was also considered in \cite{Zhao-Chen-Zhang.18} by considering  the  SCAD penalty for the    expectile process. They propose an algorithm that converges,  with probability converging to one as $n \rightarrow \infty$, to the oracle estimator after several iterations.  
Always for $p$ depending on $n$, especially when $p>n$ and for $(\varepsilon_i)$ sous-Gaussian errors,  \cite{Gu-Zou.16} penalize the expectile process with LASSO or nonconvex penalties. They find the   convergence rate of the penalized estimator, propose an algorithm   for  finding this estimator and implement  the algorithm    in the R language in package \textit{SALES}.\\
The paper of \cite{Spiegel-Sobotka-Kneib.17}    introduces several approaches depending on selection criteria and shrinkage methods to perform model selection in semiparametric expectile regression. \\

Give some general notations.  For a vector $\textbf{v}$, we denote its transpose by $\textbf{v}^t$, by $\|\textbf{v}\|_1$, $\|\textbf{v}\|_2$ and $\|\textbf{v}\|_\infty$ the $L_1$, $L_2$, $L_\infty$ norms, respectively.   The number $p$ of the explanatory variables and $p_0$ of the sugnificant variables can depend on $n$, but for convenience, we do not write the subscript $n$. Throughout the paper, $C$ denotes a positive generic constant not dependent on $n$, which may take a different value in different formula. \\
In order to study the properties of the adapted LASSO expectile estimator $\widehat  \eb_n$, we introduce the following function, using the same notations in \cite{Liao-Park-Choi.18}: 
\begin{align*}
g_\tau(x) &\equiv \rho'_\tau(x-t) |_{t=0}=2 \tau x \e1_{x \geq 0}+2(1-\tau)x \e1_{x<0}, \\
h_\tau(x)& \equiv \rho''_\tau(x-t)_{t=0}=2 \tau \e1_{x \geq 0}+2(1-\tau) \e1_{x<0}.
\end{align*}
\hh The paper is organized as follow. In Section \ref{casec<1} we study the asymptotic behavior of the adaptive LASSO expectile estimator when $p=O(n^c)$, with $0 \leq c <1$. We obtain the convergence rate of the $\widehat{\eb}_n$ and the oracle properties. A similar study is realized in Section \ref{casec>1}, when $p \geq n$. In Section \ref{Simus}, a  simulation study and an application to real data are presented. All proof are relegated in Section \ref{Proofs}.

\section{Case $c <1$, parameter number  less than the sample size}
\label{casec<1}
In this section we study the asymptotic behavior of the adaptive LASSO expectile estimator when the number $p$ of model parameter is $p=O(n^c)$, with $0 \leq c <1$. If $c=0$, that is $p$ fixed, then we get the particular case studied by \cite{Liao-Park-Choi.18}. \\
\hh An additional assumption  to   (A2)  on the design is requested:\\
\textbf{(A3)}   $p^{1/2} n^{-1/2} \max_{1 \leqslant i \leqslant n}\|\eX_i\|_2    {\underset{n \rightarrow \infty}{\longrightarrow}} 0$.\\
 
Assumption (A3) is common in works that consider $p$ dependent on $n$ when $p<n$, see for example \cite{Ciuperca-18}. Always for a linear model with the number of   parameters to order $n^c$, with $0 \leq c <1$,  \cite{Zou-Zhang-09} consider a stronger assumption for design: $\lim_{n \rightarrow \infty} n^{-1} \max_{1 \leqslant i \leqslant n} \| \eX_i\|^2_2=0$. \\
 Because $p< n$, the regression parameters are identifiable and we can calculated the expectile estimator:
\[
\widetilde   \eb_n \equiv \argmin_{\eb \in \R} Q_n(\eb),   
\]
the components of $\widetilde   \eb_n$ being $\widetilde   \eb_n=\big(\widetilde{\beta}_{n,1}, \cdots , \widetilde{\beta}_{n,p} \big)$. This estimator will intervene in the adaptive weight of the  penalty, $\widehat \omega_{n,j}  =|\widetilde \beta_{n,j}|^{-\gamma}$, for $j =1, \cdots , p$, conditions on the  constant $\gamma >0$ will be specified in Theorem \ref{Theorem 2SPL}.
By the following theorem  we obtain the convergence rate of the expectile and adaptive LASSO expectile estimators. We obtain that the convergence rate depends on the size  $p$ of the vector $\eb$.
\begin{theorem}
\label{th_vconv}
Under assumptions  (A1)-(A3) we have:\\
(i) $\| \widetilde{\eb}_n-\ebo \|_2=O_{\PP}\pth{\sqrt{\frac{p}{n}}}$.\\
(ii) If the  tuning parameter sequence $(\lambda_n)_{n \in \N}$ satisfies   $ p_0^{1/2} n^{(1-c)/2} \lambda_n \rightarrow 0$, as $n \rightarrow \infty$, then,   $\| \widehat{\eb}_n- \ebo\|_2=O_{\PP}\pth{\sqrt{\frac{p}{n}}}$.
\end{theorem}
Theorem \ref{th_vconv} provides that the expectile and adaptive LASSO expectile estimators  have the same convergence rate. Concerning the adaptive LASSO expectile estimator, the same convergence rate has been obtained for other adaptive LASSO estimators: by the likelihood method   for a generalized linear model when $p<n$  in \cite{Wang-Wang.14}, by  the least squares approximation method in \cite{Leng-Li.10}.\\
By the following theorem, considering a supplementary condition on $\lambda_n$, $c$ and $\gamma$, in addition to that considered for the convergence rate in Theorem \ref{th_vconv}, we show that the adaptive LASSO expectile estimator $\widehat{\eb}_n$ satisfies the  oracle properties. If $\tau=0.5$, that is, for the adaptive LASSO LS-estimator,   the variance of the normal limit law is the variance of $\varepsilon$. In fact, we obtained for $\tau=0.5$, the same asymptotic normality as in \cite{Zou-Zhang-09}. 
\begin{theorem}
\label{Theorem 2SPL} Suppose that assumptions  (A1)-(A3) hold and that  the tuning parameter satisfies   $\lambda_n n^{(1-c)(1+\gamma)/2} \rightarrow \infty$, $ p_0^{1/2} n^{(1-c)/2} \lambda_n \rightarrow 0$, as $n \rightarrow \infty$.  Then:\\
(i) $\PP \big[\widehat{\cal A}_n={\cal A}\big]\rightarrow 1$,  for $n \rightarrow \infty $.\\
(ii)  For any  vector $\eu$ of size $p_0$ such that $\| \eu\|_2=1$,  we have:  $n^{1/2} (\eu^t \eU^{-1}_{n,{\cal A}} \eu)^{-1/2} \eu^t ( \widehat{\eb}_n - \eb^0)_{\cal{A}}  \overset{\cal L} {\underset{n \rightarrow \infty}{\longrightarrow}} {\cal N}\big(0, \frac{\Var[g(\varepsilon)]}{\eE^2[h(\varepsilon)]} \big)$.
\end{theorem}
The convergence rate of $ \widehat{\eb}_n$ does not depend on the power $\gamma$, but otherwise, for holding the oracle properties, the  choice of $\gamma$ is very important. 
Concerning the suppositions and results stated in Theorem \ref{Theorem 2SPL}, let's make some remarks on the regularization parameter $\lambda_n$, the constants $c$, $\gamma$ and the sizes $p_0$, $p$.
\begin{remark}
1) If $p_0=O(p)$, then for that  $\lambda_n n^{(1-c)(1+\gamma)/2} \rightarrow \infty$ and $\lambda_n p_0^{1/2} n^{(1-c)/2} \rightarrow 0$ occur, we must  choose the constant $\gamma$ and sequence $(\lambda_n)$ such that: $\gamma > c/(1-c)$ and $n^{-1/2}\lambda_n \rightarrow 0$, as $n \rightarrow \infty$.\\
2) If $p_0=O(1)$, we must choose the constant $\gamma >0$ and the tuning parameter such that  $n^{(1-c)/2}\lambda_n \rightarrow 0$, as $n \rightarrow \infty$.\\
3) For that $\lambda_n p_0^{1/2} n^{(1-c)/2} \rightarrow 0$ holds, it is necessary that $\lambda_n \rightarrow 0$, as $n \rightarrow \infty$.\\
4) If $c=0$ then the conditions on $\lambda_n$ become: $n^{1/2}\lambda_n \rightarrow 0$ and $ n^{(1+\gamma)/2} \lambda_n \rightarrow \infty$, as $n \rightarrow \infty$, conditions  considered by \cite{Liao-Park-Choi.18} for a linear model with $p$ fixed. We also find the same variance of Gaussian distribution in   \cite{Liao-Park-Choi.18}.
\end{remark}

\section{Case $c \geq 1$, parameter number  greater than the sample size}
\label{casec>1}
In this section, after we propose an adaptive weight, we study the asymptotic behavior of the estimator $\widehat{\eb}_n$ when the number of regressors exceeds the number of observations. \\
\hh Since we consider now that $p \geq n$, instead of assumption (A3), we consider: \\
\textbf{(A4)}    There exists a constant $M>0$ such that $\max_{1 \leqslant i \leqslant n}\| \eX_i \|_{\infty} < M$.\\

\noindent The same assumption (A4) was considered for a generalized linear model when $p >n$ in \cite{Wang-Wang.14}  where  the adaptive LASSO likelihood method is proposed.\\
The asymptotic properties of the adaptive LASSO LS-estimators in a linear model,  when the number of the explanatory variables is greater than $n$, have been studied by \cite{Huang-Ma-Zhang.08}. They show that if a reasonable initial estimator is available, estimator that enters in the adaptive weight of the penalty,  then their  adaptive LASSO LS-estimator satisfies the oracle properties.\\
Recall that the expectile estimator is not consistent when $p>n$ and then it can't be used in the weights $\widehat{\omega}_{n,j}$ (see for example \cite{Huang-Ma-Zhang.08}). Then, when $p>n$, we propose in this section to take  as adaptive weight $\widehat{\omega}_{n,j}=\min(|\overset{\vee}{\beta}_{n,j}|^{-\gamma} , n^{1/2})$, with $\overset{\vee}{\beta}_{n,j}$ an estimator of $\eb^0_j$ consistent with $a_n \rightarrow 0$ the convergence rate: $\| \overset{\vee}{\eb}_n -\ebo \|_2=O_{\PP}(a_n)$. If the estimator $\overset{\vee}{\beta}_{n,j}$ take the value 0, then we consider $n^{1/2}$ as adaptive weight.  
An example of such estimator is the LASSO expectile estimator, proposed by  \cite{Gu-Zou.16}, defined as: 
\[ \argmin_{\eb \in \R^p}\big(n^{-1} \sum^n_{i=1} \rho_\tau(Y_i - \eX_i^t \eb)+\nu_n \| \eb\|_1 \big), 
\]
with the deterministic sequence $\nu_n \in (0, \infty)$, $\nu_n \rightarrow 0$ as $n \rightarrow \infty$.  If $\varepsilon_i$ is sub-Gaussian and $\eE[g(\varepsilon)]=0$, under our assumptions (A2), (A4), if $\kappa =\inf_{\textbf{d} \in {\cal C}} \frac{\|\eeX \textbf{d} \|^2_2}{\|n\textbf{d}\|^2_2} \in (0, \infty)$, with $\eeX$ the matrix $n \times p $ of design and the set ${\cal C} \equiv \{\textbf{d} \in \R^p ; \|\textbf{d}_{{\cal A}^c}\|_1 \leq 3 \|\textbf{d}_{{\cal A}}\|_1 \neq 0\}$, then $\|\overset{\vee}{{\eb}}_n -\ebo \|_2=O_{\PP}\big(p_0^{1/2} \nu_n  \big)$. Thus, the sequence $(a_n)$ is in this case $a_n=p_0^{1/2} \nu_n  $ (see Theorem 1 of  \cite{Gu-Zou.16}).\\

The form of the random process $R_n(\eb)$ of (\ref{Rnb}) and the adaptive LASSO expectile estimator $\widehat{\eb}_n$ of (\ref{hbetan}) remain the same, only the adaptive weight changes. It would be desirable for $\widehat{\eb}_n$ to satisfy the oracle proprties. 
For the sparsity  property of   $\widehat{\eb}_n$ its  convergence in $L_1$  norm   is required.  In the following theorem, $(b_n)_{n \in \N}$ is a deterministic sequence converging to 0 as $n \rightarrow \infty$.

\begin{theorem}
\label{th_vconvn}
Under assumptions  (A1), (A2), (A4), the  tuning parameter $(\lambda_n)_{n \in \N}$ and sequence $(b_n)$  satisfying  $\lambda_n p_0^{1/2} b^{-1}_n \rightarrow 0$, as $n \rightarrow \infty$, we have,   $\| \widehat{\eb}_n- \ebo\|_1=O_{\PP}\pth{b_n}$.
\end{theorem}
The result of Theorem \ref{th_vconvn} indicates that the convergence rate of the adaptive LASSO expectile  estimator $\widehat{\eb}_n$ depends on the chosen sequence  $(\lambda_n)_{n \in \N}$. On the other hand, the convergence rate of $\widehat{\eb}_n$ don't depend on the convergence rate $(a_n)$ of the   estimator $\overset{\vee}{{\eb}}_n$. The only thing that matters (see relation (\ref{PPn}) of the proof in Section \ref{Proofs}) is that  $\overset{\vee}{{\eb}}_n$ converges in probability to   $\ebo$.\\
The result of Theorem \ref{th_vconvn} now allows us to state oracle properties.  Still in the case $p>n$, but for a quantile model, \cite{Zheng-Gallagher-Kulasekera-13} obtains that the convergence rate, in $L_2$ norm, of the adaptive LASSO quantile estimator is $(p_0 /n)^{-1/2}$ and that it also satisfies the oracle properties. The sparsity of the  adaptive LASSO quantile estimator will be shown in our Section \ref{Simus} by  simulations, where we obtain that compared to the adaptive   LASSO expectile estimator, it would be necessary to have a larger number $n$ of the observations when the  model errors have an asymmetric distribution. From where, a supplementary interest in considering the expectile method instead of the quantile.
\begin{theorem}
\label{Theorem 2SPLn} Suppose that assumptions   (A1), (A2), (A4) hold,   the tuning parameter $(\lambda_n)$ and sequence $(b_n)$ satisfy   $\lambda_n p_0^{1/2} b^{-1}_n \rightarrow 0$, $\lambda_n b_n^{-1}\min\big(n^{1/2} , a_n^{-\gamma}\big) \rightarrow \infty$, as $n \rightarrow \infty$.  Then:\\
(i)  $\PP \big[\widehat{\cal A}_n = {\cal A}\big]\rightarrow 1$,  for $n \rightarrow \infty $.\\
(ii)   For any  vector $\eu$ of size $p_0$ such that $\| \eu\|_1=1$,   we have  $n^{1/2} (\eu^t \eU^{-1}_{n,{\cal A}} \eu)^{-1/2} \eu^t ( \widehat{\eb}_n - \eb^0)_{\cal{A}}  \overset{\cal L} {\underset{n \rightarrow \infty}{\longrightarrow}} {\cal N}\big(0, \frac{\Var[g(\varepsilon)]}{\eE^2[h(\varepsilon)]} \big)$.
\end{theorem}
Hence, even if the parameter number of the model is larger than the observation number, the variance of the    normal limit distribution is the same as that obtained when $p<n$. As for the case $p < n$, studied in Section \ref{casec<1}, the convergence rate of the adaptive LASSO espectile estimator don't depend on the power $\gamma$  in the adaptive weight. However, $\gamma$ intervenes in the imposed conditions so that the oracle properties are satisfied. If $\tau=0.5$, that is for the adaptive LASSO LS-estimator, we obtained the same asymptotic normal distribution as that given by \cite{Huang-Ma-Zhang.08}) for their adaptive LASSO LS-estimator.\\
Regarding the tuning parameter sequence we make the following remark, useful for simulations and applications on real data.
\begin{remark}
The supposition $\lambda_n p_0^{1/2} b^{-1}_n \rightarrow 0$ made in the Theorems \ref{th_vconvn} and \ref{Theorem 2SPLn}, implies  that the  tunning sequence $\lambda_n \rightarrow 0$, as $n  \rightarrow \infty$.
\end{remark}

\section{Numerical study}
\label{Simus}
In this section we first perform a numerical simulation study to illustrate our theoretical results  on the adaptive LASSO expectile estimation and  compare it with estimation obtained by the adaptive LASSO quantile method. Afterwards, an application on real data is presented.\\
We use the following R language packages: package \textit{SALES} (for $p>n$) with function \textit{ernet} for the expectile regression and for quantile regression, the package \textit{quantreg} with function \textit{rq}.\\
Given assumption  (A1), the expectile index $\tau$ is:
\begin{equation}
\label{te}
\tau=\frac{\eE\big[\varepsilon \e1_{\varepsilon<0} \big]}{\eE\big[\varepsilon(\e1_{\varepsilon<0}- \e1_{\varepsilon>0})  \big]}.
\end{equation}
 In the simulation study, the index  $\tau$ is fixed, function of the law of $\varepsilon$, such that assumption (A1) be satisfied. On basis relation (\ref{te}), we will give details in subsection \ref{appli} how an estimation for $\tau$  can be calculated in pratical applications.\\
Taking into account the suppositions imposed on the tuning parameter in Theorems \ref{Theorem 2SPL} and  \ref{Theorem 2SPLn}, we consider   $\lambda_n=n^{-2/5}$ for the  expectile framework.   For the power  $\gamma$ in the adaptive weights of the penalty,  several values will be considered and a variation on a grid of values will be realized to choose the values of $\gamma$ which give the best results in terms of the significant variable identification. 
 For the quantile method, the tuning parameter is $n^{2/5}$ and  the weight in the penalty have the power $1.225$ (see \cite{Ciuperca-16}).
 \subsection{\textbf{Simulation study: fixed $p_0$ case}}
  In this subsection, we will study the numerical behavior of the adaptive LASSO expectile method   and we will compare it with the simulation results   obtained by the adaptive LASSO quantile method.
  \begin{table} 
\caption{\footnotesize  Sparsity study   for expectile   method with adaptive LASSO penalty, $\varepsilon \sim {\cal N}(0,1)$, two values for $\gamma$. }
\begin{center}
{\tiny
\begin{tabular}{|c|c|cc|cc|}\hline 
  $n$ & $p$ &   \multicolumn{2}{c|}{$Card({\cal A} \cap \widehat{{\cal A}}_n)$} & \multicolumn{2}{c|}{$Card(\widehat{{\cal A}}_n \setminus {\cal A})$} \\ \cline{3-6}
   & & $\gamma=1/8$ &$\gamma=2$ &   $\gamma=1/8$ & $\gamma=2$   \\ \hline
  50 & 100 & 5.99&5.58   & 3.8&0.004  \\ 
   & 400 & 5.87&5.12   & 8.7&0.04  \\
    & 600 & 5.71&4.98   & 10&0.07  \\ \hline
     100 & 100 & 6&5.99   & 2.48&0   \\ 
   & 400 & 6&5.96   & 6.79&0   \\
    & 600 & 6&5.97   & 8.8&0  \\ \hline
   \end{tabular}
}
 \end{center}
\label{Tabl0} 
\end{table} 
 For model (\ref{eq1}), we consider $p_0=6$ and  ${\cal A}=\{1, \cdots , 6 \}$.  In the all simulations of this subsection we take, $\beta^0_1=1$, $\beta^0_2=4$, $\beta^0_3=-3$, $\beta^0_4=5$, $\beta^0_5=6$, $\beta^0_6=-1$, while $n$ and $p$ are varied. The values of $p$ they can be less than $ n $ but also higher. \\
For the errors $\varepsilon_i$, three distributions are considered: ${\cal N}(0,1)$ which is  symmetrical, ${\cal E}xp(-2.5)$ and ${\cal N}(0,4 \cdot 10^{-2})+\chi^2(1)$, the last two being asymmetrical. The exponential law  ${\cal E}xp(-2.5)$ has the density function  $\exp(-(x+2.5))\e1_{x >-2.5}$. For each value of  $n$, $p$ and  distribution of $\varepsilon$, 1000   Monte Carlo replications are realized for two possible values for $\gamma$.    In  Tables \ref{Tabl0} and  \ref{Tabl1} we give the average of the 1000 Monte Carlo replications for the cardinalities (number of the true non-zeros estimated as non-zero) $Card({\cal A} \cap \widehat{{\cal A}}_n)$ and $Card(\widehat{{\cal A}}_n \setminus {\cal A})$ (number of the false non-zero) by the   expectile (ES) and quantile (Q)  penalized methods, each with LASSO adaptive penalty. For a perfect method, we should have: $Card({\cal A} \cap \widehat{{\cal A}}_n)=p_0=6$ and $Card(\widehat{{\cal A}}_n \setminus {\cal A})=0$. In Table \ref{Tabl0}, for standard Gaussian errors, the values considered for $\gamma$ are $1/8$ and $2$. Since for $\gamma=1/8$ the number of false non-zeros, which in addition increases with $n$, is much larger than for $\gamma=2$ and for $\gamma=2$ the number of the true non-zeros decreases with $n$, these values will be dropped, other two values will be considered in Table   \ref{Tabl1}. In Table \ref{Tabl1},  taking $\gamma \in \{5/8, 1\}$, for the adaptive LASSO expectile method, all significant variables are detected when $ n $ is  large enough, the number $ p $ of variables not coming into play.   
\begin{table} 
\caption{\footnotesize  Sparsity study   for expectile (ES) and quantile (Q) methods with adaptive LASSO penalties. }
\begin{center}
{\tiny
\begin{tabular}{|c|c|c|ccc|ccc|}\cline{1-9}
$\varepsilon$ & $n$ & $p$ &   \multicolumn{3}{c|}{$Card({\cal A} \cap \widehat{{\cal A}}_n)$} & \multicolumn{3}{c|}{$Card(\widehat{{\cal A}}_n \setminus {\cal A})$} \\ \cline{4-9}
& & & \multicolumn{2}{c}{ ES} & Q &\multicolumn{2}{c}{ ES} & Q \\  \cline{4-5} \cline{7-8}
 & & & $\gamma=5/8$ &$\gamma=1$ & & $\gamma=5/8$ & $\gamma=1$ & \\ \cline{1-9}
 $ {\cal N} (0,4 \cdot 10^{-2})+\chi^2(1)$ & 50 & 100 &   5.25 & 5&5.11& 0.06 &0.03 & 0.16\\ 
  & & 400 & 4.98 & 4.67  & 4.45 &  0.34 & 0.12 & 0.59 \\ 
    & & 600 & 4.80  &  4.57 & 4.32 &  0.54 & 0.19& 0.67 \\ \cline{2-9}
    & 100 & 100 &  5.96  & 5.85 & 5.98&   0   &0 & 0 \\ 
  & & 400 & 5.94  & 5.85 & 5.98 &  0.003 & 0& 0.002 \\ 
    & & 600 &  5.94 & 5.84& 5.96 &  0.006  &0 &  0.001  \\  \cline{2-9}
      & 200 & 100 &   6&6 &6&  0 &0 & 0 \\ 
  & & 400 & 6 &6 &6 & 0& 0& 0  \\ 
    & & 600 &  6 &6 &6& 0& 0 & 0  \\  \cline{1-9}
     $ {\cal E}xp(-2.5)$ & 50 & 100 & 5.8 &5.66  & 3.18 &  0.33&0.06 & 0.57 \\ 
  & & 400 &  5.58&5.41 & 2.9 & 0.93& 0.27 & 1.1 \\ 
    & & 600 & 5.37  & 5.11 &2.7 &  1.3&0.46 & 1.25 \\  \cline{2-9}
    & 100 & 100 & 6  & 5.99&4.9  & 0.06 &0.003 & 0.09  \\ 
  & & 400 &  5.99& 5.99 & 4.67 &   0.12&0.008 & 0.12  \\ 
    & & 600 & 6& 5.99 & 4.68 &  0.14&0.02 & 0.10   \\  \cline{2-9}
      & 200 & 100 &   6& 6 & 5.77 &  0.007&0.002  &  0.04 \\ 
  & & 400 & 6& 6 &  5.79 &  0.01&  0 &  0.02 \\ 
    & & 600 &  6& 6 & 5.79 & 0.019& 0  &  0.020  \\  \cline{1-9}
         $ {\cal N}(0,1)$ & 50 & 100 & 5.95&5.88  & 5.89 & 0.21&0.02 & 0.44 \\ 
  & & 400 &  5.74 & 5.58 &5.51 &  0.72& 0.28 & 0.72  \\ 
    & & 600 &  5.61&5.36  &  5.33  &  1.2& 0.40 & 0.92 \\  \cline{2-9}
    & 100 & 100 &  6 & 6 &6  &  0.01& 0&  0.14  \\ 
  & & 400 &  6& 6  & 6 & 0.04&0.006  & 0.17   \\ 
    & & 600 &  5.99& 6   & 5.99  &  0.04& 0  &   0.15  \\  \cline{2-9}
      & 200 & 100 &   6&  6 &  6 & 0.001 & 0 &  0.03  \\ 
  & & 400 &  6& 6 &  6  &   0.001 & 0 &0.06   \\ 
    & & 600 &  6& 6 & 6 &  0  &  0 &0.05  \\  \cline{1-9}
        \end{tabular}
}
 \end{center}
\label{Tabl1} 
\end{table}  
We observe that the penalized expectile method better detects non-zero parameters compared to the penalized quantile method, especially for exponential errors. For small values of $ n $, the expectile method makes less false detections of non-significant variables as significant variables. Concerning the two values considered for $\gamma$, when $\gamma = 5/8$, there are a little more true non-zeros detected, while, when $\gamma=1$, there are fewer false non-zeros. This trend will be also confirmed by the following numerical studies.
 \subsection{\textbf{Simulation study: case when $p_0$ varies with $n$}}
 In this subsection, we always compare  expectile and quantile penalized methods, but when the values considered for $ p $    are larger than $ n $. Moreover, the number of non-zero parameters can increase as $n$ increases. 
 \begin{table} 
\caption{\footnotesize Study of the sparsity evolution and of the estimation  accuracy for the expectile (ES) and quantile framework, with $p$ and $p_0$ depending on  $n$: $p=4n$, $p_0=2[n^{1/2}]$. }
\begin{center}
{\tiny
\begin{tabular}{|c|c|ccc|ccc|ccc|ccc|}\cline{1-14}
$\varepsilon$ &  $n$  &   \multicolumn{3}{c|}{$100  p_0^{-1}   Card({\cal A} \cap \widehat{{\cal A}}_n)$} & \multicolumn{3}{c|}{$100   (n-p_0)^{-1}  Card(\widehat{{\cal A}}_n \setminus {\cal A})$}  & \multicolumn{3}{c|}{$mean(|\widehat{\eb}_n-\ebo|)$}  & \multicolumn{3}{c|}{$mean(|(\widehat{\eb}_n-\ebo\big)_{\cal A}|)$}\\ \cline{3-14}
 & & \multicolumn{2}{c}{ ES} & Q &\multicolumn{2}{c}{ ES} & Q & \multicolumn{2}{c}{ ES} & Q &\multicolumn{2}{c}{ ES} & Q\\  \cline{3-4} \cline{6-7} \cline{9-10} \cline{12-13}
  & & $\gamma=5/8$ &$\gamma=1$ & & $\gamma=5/8$ & $\gamma=1$ & & $\gamma=5/8$ &$\gamma=1$ & & $\gamma=5/8$ & $\gamma=1$ & \\ \hline \hline
$ {\cal N}(0,1)$ &  75 & 98 &97.4 &98 & 3.3 & 2.7&3.6 & 0.01 &0.01 &0.01 & 0.22& 0.23& 0.22\\ 
 & 100 & 99.4 & 99.2 & 99.3 & 1.1 & 0.46 &1.5 & 0.007 & 0.006 &0.007  & 0.13 &0.13 &0.13 \\ 
& 200 & 100 &100 &100 & 0.02 & 0.008 &0.07  & 0.002 & 0.002&0.002 & 0.07 & 0.07&0.07\\ 
 &400 & 100 &100 &100 & 0.01 &0.008 &0.002 & 0.001&0.001 &0.001 & 0.05 &0.04 &0.05 \\ \hline
 $\varepsilon \sim  {\cal E}xp(-2.5)$  &  75 & 97 &97 &43 & 4.2 &2.21 &10.9 & 0.01 & 0.01&0.34 &0.26 &0.23 &5.2\\ 
   & 100 & 99 &98.6 &49 & 1.6 & 0.78&11.6 & 0.008 &0.008 &0.35  & 0.16 &0.15 &5.7 \\ 
   & 200 & 100 & 100&77 & 0.05& 0.01&5  & 0.002 & 0.002&0.13  & 0.08 & 0.08&3.1\\ 
    &400 & 100 & 100 &98 & 0.04 & 0.02&0.37 & 0.001&0.001 &0.007 & 0.05 & 0.05&0.28\\ \cline{1-14}
    ${\cal N} (0,4 \cdot 10^{-2})+\chi^2(1)$  &  75 & 94 &93.3 &70 & 3.7 &2.55 &12.4 & 0.02 &0.02 &0.18 &0.40& 0.37& 2.58\\ 
     & 100 & 97&  96 &80 & 1.5 & 0.81&9.8 & 0.01 & 0.01&0.13  & 0.26 & 0.22&2.07 \\
      & 200 & 99.8 & 99.6&99 & 0.02&0.01 &0.26  & 0.003 & 0.002&0.003  & 0.1 & 0.08&0.1\\ 
     &400 & 100 &100 &100 & 0.02 &0.002& 0 & 0.001& 0.001&0.001 & 0.05& 0.04& 0.02\\ \cline{1-14}
 \end{tabular}
}
 \end{center}
\label{Tabl2} 
\end{table}  
In Table \ref{Tabl2} we take $p=4n$, $p_0=2[n^{1/2}]$, with $[x]$ the entire part of $x$, the power  $\gamma \in \{5/8, 1\}$ and $\varepsilon \sim  {\cal E}xp(-2.5)$. The true value of the non-null parameter vector is  $\eb^0_{\cal A}=(1, \cdots , p_0)$. We assess model selection by calculating  the percentage ($100   p_0^{-1}  Card({\cal A} \cap \widehat{{\cal A}}_n)$) of the non-zero parameters with a non-zero estimation and the percentage of false significant variables ($100 (n-p_0)^{-1}  Card(\widehat{{\cal A}}_n \setminus {\cal A})$), by the two estimation methods. We also give the accuracy of the complete estimation vectors ($mean(|\widehat{\eb}_n-\ebo|)$) and of the estimations of non-zero parameters ($mean(|(\widehat{\eb}_n-\ebo\big)_{\cal A}|)$) (average absolute estimation error) obtained on  1000   Monte Carlo replications. More precisely, if $M$ is the Monte Carlo replication number and $\widehat{\eb}^{(m)}_{n,j}$ is the estimation of $\beta^0_j$ obtained for the Monte Carlo replication with the number $m$, then, $mean(|\widehat{\eb}_n-\ebo|)=(Mp)^{-1} \sum^{M}_{m=1} \sum^p_{j=1}| \widehat{\eb}^{(m)}_{n,j} - \beta^0_j | $. Similarly we calculate $mean(|(\widehat{\eb}_n-\ebo)_{\cal A}|)=(Mp)^{-1} \sum^{M}_{m=1} \sum^{p_0}_{j=1} | \widehat{\eb}^{(m)}_{n,j} - \beta^0_j | $. For the Gaussian standardized  errors, the results are similar by the two estimation methods. The results remain very good when the errors have an exponential law   or a mixing between a Gaussian with a $\chi^2(1)$ law for the proposed    adaptive LASSO expectile method, even for small values for $n$,  the results being more accurate for $\gamma=1$ than for $\gamma =5/8$. Furthermore, the adaptive LASSO quantile method, provides less accurate estimations, even when it correctly detects significant variables (for values of $ n $ large), which can pose problems of application in practice.
\begin{table} 
\caption{\footnotesize Study of the sparsity evolution and of the estimation  accuracy for the expectile (ES) and quantile framework, with $p$ and $p_0$ depending an $n$: $p=[n\log(n)]$, $p_0=2[n^{1/4}]$, $\varepsilon \sim  {\cal E}xp(-2.5)$. }
\begin{center}
{\tiny
\begin{tabular}{|c|c|ccc|ccc|ccc|ccc|}\hline 
$\eb^0_{\cal A}$ &  $n$  &   \multicolumn{3}{c|}{$100   p_0^{-1}   Card({\cal A} \cap \widehat{{\cal A}}_n)$} & \multicolumn{3}{c|}{$100   (n-p_0)^{-1}  Card(\widehat{{\cal A}}_n \setminus {\cal A})$}  & \multicolumn{3}{c|}{$mean(|\widehat{\eb}_n-\ebo|)$}  & \multicolumn{3}{c|}{$mean(|(\widehat{\eb}_n-\ebo\big)_{\cal A}|)$}\\ \cline{3-14}
& & \multicolumn{2}{c}{ ES} & Q &\multicolumn{2}{c}{ ES} & Q & \multicolumn{2}{c}{ ES} & Q &\multicolumn{2}{c}{ ES} & Q\\   \cline{3-4} \cline{6-7} \cline{9-10} \cline{12-13}
& & $\gamma=5/8$ &$\gamma=1$ & & $\gamma=5/8$ & $\gamma=1$ & & $\gamma=5/8$ &$\gamma=1$ & & $\gamma=5/8$ & $\gamma=1$ & \\ \hline \hline
$(1, 2, \cdots, p_0)$ & 75 & 99.8 & 99&76 & 0.27 &0.02 &0.22 &  0.002 & 0.002&0.008 & 0.19& 0.19&0.67 \\ 
&  100 & 100 &99.9& 87 & 0.12 &0.008 &0.11 & 0.002& 0.002 & 0.007 & 0.15&0.14 & 0.53 \\ 
 & 200 & 100 & 100&98 & 0&0  & 0  & 0.0006&  0.0005& 0.001 & 0.10& 0.09& 0.24 \\ 
 & 400 & 100&100 & 100 & 0& 0 & 0  & 0.0002& 0.0005 & 0.0004 & 0.06&0.09 & 0.14 \\ \hline
  $(1, \cdots, 1)$ & 75 & 98.4 &97.9&30 & 0.32 &0.03 &0.25 &  0.004& 0.004& 0.01 & 0.31& 0.35&0.83 \\ 
  &  100 & 99.8&99.4 & 36 & 0.26& 0.03& 0.17 & 0.003& 0.003 & 0.01 & 0.25& 0.28& 0.81 \\ 
   & 200 & 100& 100& 82 & 0.02&0 & 0.005  & 0.0009& 0.001 & 0.003 & 0.16& 0.18& 0.51 \\ 
    & 400 & 100&100 & 100 & 0 & 0& 0  & 0.0003 & 0.0003& 0.007 & 0.11&0.11 & 0.23 \\ \cline{1-14}
  \end{tabular}
}
 \end{center}
\label{Tabl3} 
\end{table}  
 In Table  \ref{Tabl3}, taking $p=[n \log(n)]$, the value of  $p$ is increased compared to that considered in Table \ref{Tabl2}. Furthermore, the sparsity of the model is more accentuated by considering $p_0=2[n^{1/4}]$. Two values for $\eb^0_{\cal A}$ are considered: $(1, 2, \cdots, p_0)$ and $\textbf{1}_{p_0}=(1,\cdots , 1)$ while  for the model errors, only the exponential distribution   $\varepsilon \sim  {\cal E}xp(-2.5)$ is made. For both values of $\ebo$, the expectile method with  adaptive LASSO penalty  gives very good results for identifying of null and non-null parameters, while the quantile method identifies all significant variables only when $ n $ is large (greater than 200). \\
 Comparing Tables \ref{Tabl2} and \ref{Tabl3}, for $n$ fixed and exponential law,  we deduce that the penalized expectile estimation  quality of the model does not vary for two different $ p $. Furthermore, the quality is better if $ p_0 $ decreases and when  $\gamma =1$.
 \subsection{\textbf{Sparsity study function of $\gamma$}}
 \begin{figure}[h!]
\includegraphics[width=16.5cm,height=6cm,angle=0]{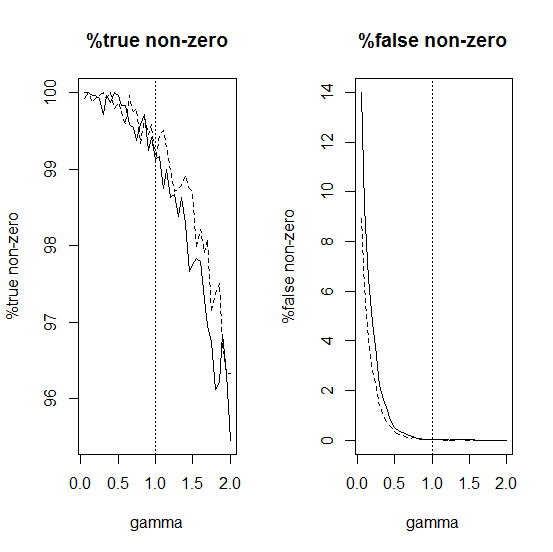}
\caption{ \it   Percentage of true (to the left) and false non-zero, for  $\varepsilon \sim  {\cal E}xp(-2.5)$, $n=75$. Two value for the number of parameters: $p=100$ (doted line) and $p=200$ (solid line).}
\label{Figure n75_exp}
\end{figure}
\begin{figure}[h!]
\includegraphics[width=16.5cm,height=6cm,angle=0]{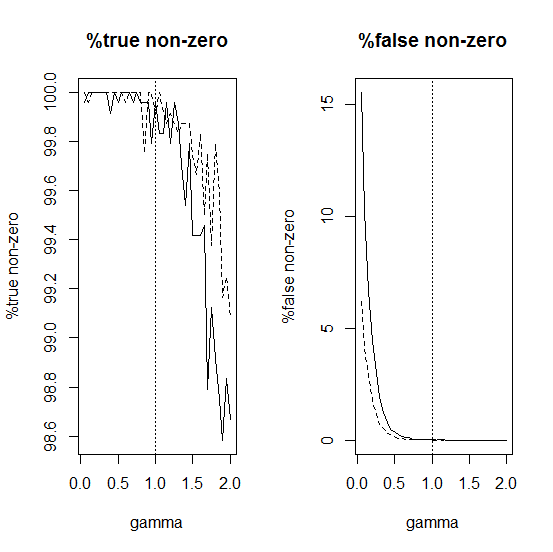}
\caption{ \it   Percentage of true and false non zero, for  $\varepsilon \sim  {\cal E}xp(-2.5)$, $n=100$. Two value for the number of parameters: $p=100$ (doted line) and $p=400$ (solid line).}
\label{Figure n100_exp}
\end{figure}
\begin{figure}[h!]
\includegraphics[width=16.5cm,height=6cm,angle=0]{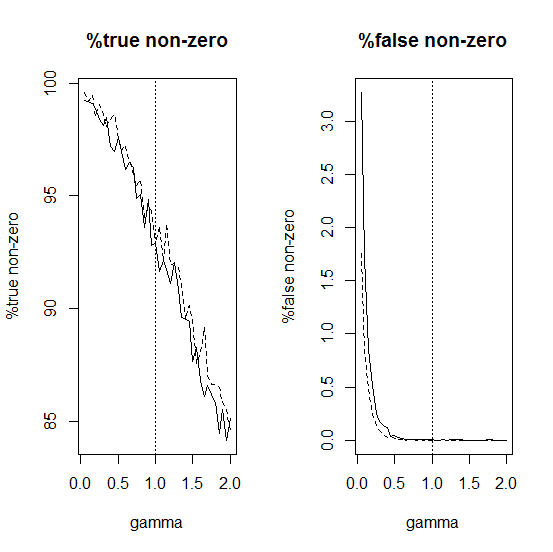}
\caption{ \it   Percentage of true and false non zero, for  $\varepsilon \sim  {\cal N} (0,4 \cdot 10^{-2})+\chi^2(1)$, $n=75$. Two value for the number of parameters: $p=100$ (doted line) and $p=200$ (solid line).}
\label{Figure n75_mixt}
\end{figure}
\begin{figure}[h!]
\includegraphics[width=16.5cm,height=6cm,angle=0]{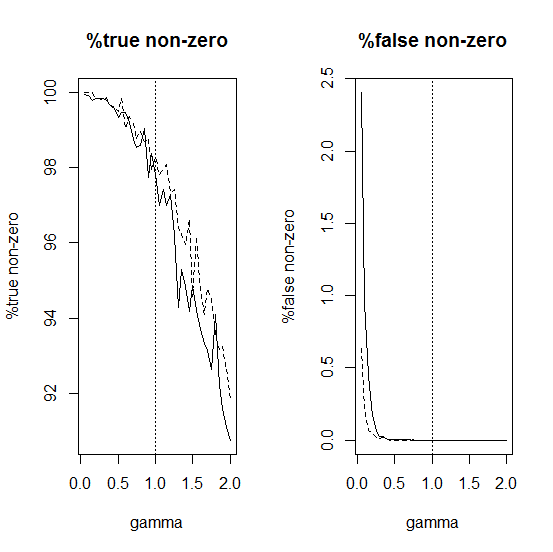}
\caption{ \it   Percentage of true and false non zero, for  $\varepsilon \sim  {\cal N} (0,4 \cdot 10^{-2})+\chi^2(1)$, $n=100$. Two value for the number of parameters: $p=100$ (doted line) and $p=400$ (solid line).}
\label{Figure n100_mixt}
\end{figure}
 In Figure \ref{Figure n75_exp} we present for a model with sample size  $n=75$, errors $\varepsilon \sim  {\cal E}xp(-2.5)$, the percentage of the true   non-zero parameters estimated by adaptive LASSO expectile method (sub-figure to the left) and the percentage of the true zeros estimated as non-zero (to right), for two values of $p$: $p=100$ (dotes line) and $p=200$ (solid line). We observe that for $\gamma \in ( 0, 1.2]$ the estimation rate of the non-zero parameters by a non-zero estimation exceeds $99 \%$, for $\gamma \geq 0.5$ the rate of false non-zero   is less than $1\%$ and this last rate decreases when $\gamma $ increases. A similar study is presented in Figure \ref{Figure n100_exp} for $n=100$, $p \in \{100, 400\}$ and we deduce the same conclusions as in Figure \ref{Figure n75_exp}. \\
In view of all this, we may deduce that for a power  $\gamma \in [0.5, 1.25]$, the error rates (of false zeros and false non-zeros) are less than $1 \%$.\\
 For errors  $\varepsilon \sim  {\cal N} (0,4  \cdot 10^{-2})+\chi^2(1)$ (see Figures \ref{Figure n75_mixt} and \ref{Figure n100_mixt}), the detection rate of the true non-zero parameters decreases faster when $\gamma$  increases, the rate of false non-zeros remaining the same as for exponential errors (Figures \ref{Figure n75_exp} and \ref{Figure n100_exp}).
 \subsection{\textbf{Application to real data}}
 \label{appli}
\begin{figure}[h!]
\includegraphics[width=16.0cm,height=6cm,angle=0]{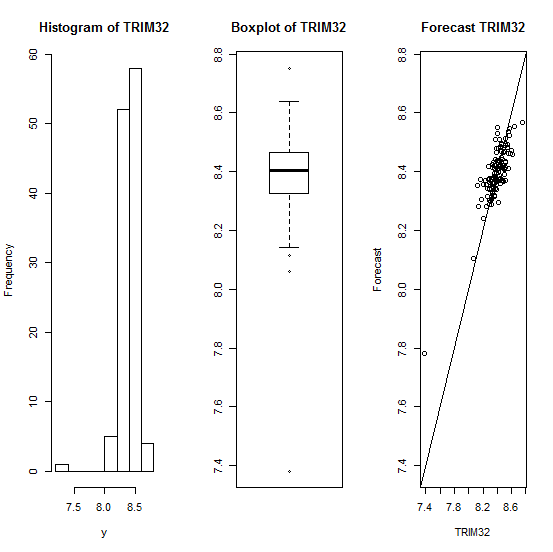}
\caption{ \it   Histogram and boxplot of TRIM32. Scater graph between forecast and the true value of  TRIM32.}
\label{Figure 1}
\end{figure}

 We use the data \textit{eyedata} of  R package  \textit{flare} which contains $n=120$ observations (rats) for the response variable  of gene \textit{TRIM32} and 200 explanatory variables, other genes probes, from the microarray experiments of mammalian-eye tissue samples in   \cite{Scheetz.06}. The objective is to find genes that are correlated with the TRIM32 gene, known to
cause Bardet–Biedl syndrome, a genetically disease of multiple organ systems including the retina.\\
In order to calculate in applications the expectile index  $ \tau $, we  standardize the values of the explained variable $\widetilde{y}_i=(y_i-\bar y_n)/(\widehat{\sigma}_y)$. Afterwards, based on relation (\ref{te}), we calculate the empirical estimation of $\tau$:
 \[
 \widehat \tau=\frac{n^{-1}\sum^n_{i=1} \widetilde{y}_i \e1_{\widetilde{y}_i <0}}{n^{-1} \big( \sum^n_{i=1} \widetilde{y}_i \e1_{\widetilde{y}_i <0}-\sum^n_{i=1} \widetilde{y}_i \e1_{\widetilde{y}_i >0}\big)}.
 \]
Then, in model (\ref{eq1}), the response variable is  $\widetilde{Y}_i=(Y_i-\overline{Y}_n)/\widehat{\sigma}_Y$, with $\widehat{\sigma}_Y$ the empirical standard deviation and $\overline{Y}_n$ empirical mean of $Y$. For this application, we get $ \widehat \tau=0.533$ and $\gamma=5/8$.\\
 We obtain, for $\gamma = 5/8$, $\lambda_n=n^{2/5}$, by the method adaptive LASSO expectile that the genes whose expressions influence gene TRIM32 are 87, 153, 180, 185 with the labels: "21092", "25141", "28680" and "28967". The obtained estimations for the  coefficients of these four explanatory variables are respectively: \textit{-0.65, 2.24, 0.28} and \textit{-0.35}. In Figure \ref{Figure 1} we illustrate the histogram and the boxplot for  response variable TRIM32. We observe that it don't have a symmetrical law.  \\
If a classical LS regression  of the TRIM32 variable in respect to  the four selected covariates is performed, we obtain a model with an adjusted $ R ^ 2 = 0.74 $ and with residual standard error = 0.50. The all four variables are  significant and the residuals have Gaussian distribution (the p-value by Shapiro test equal to 0.45). In the sub-figure of the right-hand side of Figure \ref{Figure 1}, we also present,   the forecasts beside of the true values of TRIM32. We observe that the scater graph is on the first  bisectrix.\\
By the adaptive  LASSO quantile method, no variable is selected among the 200 explanatory ones. \\
In literature works that model the same data, variable number 153, tagged "25141", has been selected as the sole regressor by the bayesian shrinkage in  \cite{Song-Liang.17} and by a globally adaptive quantile method in \cite{Zheng-Peng-He.15} for quantile index between 0.45 and 0.55. In this last paper, there are other covariates that appear to be significant for other quantile index  values. These variables are: "11711", "24565", "25141", "25367", "21092", "29045", "25439", "22140", "15863" and "6222". If we make a classic regression for these ten regressors, we obtain, with a risk of $ 0.05 $ that only the variables "25141", "21092", "29045", "15863", "6222" are significant, in a model of lower quality (adjusted $R^2=0.72$, residual standard error=0.52)  than the one with the four explanatory variables   found in the present paper.

 \section{Proofs}
 \label{Proofs}
 In this section we give  the proofs of the results presented in Sections \ref{casec<1} and \ref{casec>1}.
 \subsection{\textbf{Result proofs  in Section \ref{casec<1}}}
 \noindent {\bf Proof of Theorem \ref{th_vconv}}.\\
\textit{(i)} In order to show the  convergence rate of the expectile estimator, we show that for all $\epsilon >0$, there exists $B_\epsilon >0$, large enough when $ n $ is large, such that:
\begin{equation}
\label{eq2}
\PP \big[ \inf_{\eu \in \R^p, \; \| \eu \|_2 =1} Q_n\big(\ebo+B_\epsilon \sqrt{\frac{p}{n}} \eu\big)>Q_n(\ebo)\big] \geq 1-\epsilon .
\end{equation}
This part of proof is similar to that of Lemma 1.2 in  \cite{Zhao-Chen-Zhang.18}, for the  convergence rate of the  oracle estimator. 
Let $B>0$ be  a constant to be determined later and $\eu$ a vector in $\R^p$ with the norm $\| \eu \|_2=1$.  Let's study the difference:
\begin{align*}
&Q_n\big(\ebo+B \sqrt{\frac{p}{n}}  \eu \big) -Q_n(\ebo) = \sum^n_{i=1}\big[\rho_\tau\big(\varepsilon_i -B \sqrt{\frac{p}{n}} \eX^t_i  \eu \big)- \rho_\tau(\varepsilon_i ) \big]
\end{align*}
\begin{align*}
& \quad = \sum^n_{i=1}\big[\rho_\tau\big(\varepsilon_i -B \sqrt{\frac{p}{n}} \eX^t_i  \eu \big)- \rho_\tau(\varepsilon_i ) - \eE\big[ \rho_\tau\big(\varepsilon_i -B \sqrt{\frac{p}{n}} \eX^t_i  \eu \big)- \rho_\tau(\varepsilon_i ) \big]\big]  \\
& \qquad \qquad +\sum^n_{i=1} \eE\big[ \rho_\tau\big(\varepsilon_i -B \sqrt{\frac{p}{n}} \eX^t_i  \eu \big)- \rho_\tau(\varepsilon_i ) \big] \equiv \Delta_1+\Delta_2.
\end{align*}
We first study the term $\Delta_2$. By Taylor expansion, we have: $\eE[\rho_\tau(\varepsilon-t)-\rho_\tau(\varepsilon)]=\eE \big[ - g_\tau(\varepsilon)t +2^{-1} h_\tau(\varepsilon) t^2 \big]+o(t^2)=2^{-1}\eE \big[ h_\tau(\varepsilon) \big]t^2  +o(t^2)$. By the Cauchy-Schwarz inequality, we have that $| \eX^t_i \eu |^2 \leq \| \eX_i\|^2_{2} \| \eu \|_2^2$ and then, using assumption   (A3), we obtain that $(pn^{-1})^{1/2} \max_{1 \leqslant i \leqslant n}| \eX^t_i \eu | {\underset{n \rightarrow \infty}{\longrightarrow}} 0$.  Thus, 
\[
\Delta_2 =\frac{1}{2} \sum^n_{i=1}\bigg[ B^2\frac{p}{n} \big(\eX_i^t \eu \big)^2\eE \big[ h_\tau(\varepsilon) \big]+o\bigg( B^2\frac{p}{n} \big(\eX_i^t \eu \big)^2\eE \big[ h_\tau(\varepsilon) \big] \bigg) \bigg].
\]
On the other hand, 
\[
2\min(\tau, 1-\tau) \leq \eE \big[ h_\tau(\varepsilon) \big]=2 \tau \eE[\e1_{\varepsilon \geq 0}]+2(1-\tau)\eE[\e1_{\varepsilon<0}] \leq 2\max(\tau, 1-\tau).
\]
Then
\begin{equation}
\label{DD2}
0< \Delta_2 =B^2 \frac{p}{n} \eE \big[ h_\tau(\varepsilon) \big]  \sum^n_{i=1} \big( \eX_i^t \eu \big)^2\big(1+o(1)  \big)=O(B^2p).
\end{equation}
We are now studying the term $\Delta_1$. Let us consider the random variable:
\begin{equation}
\label{Ri}
D_i \equiv \rho_\tau \big(\varepsilon_i - B \sqrt{\frac{p}{n}} \eX_i^t \eu \big) - \rho_\tau(\varepsilon_i) +B \sqrt{\frac{p}{n}} g_\tau(\varepsilon_i) \eX_i^t \eu.
\end{equation}
Then, we can write $\Delta_1$ as:
\[
\Delta_1 = \sum^n_{i=1} \big[ - B \sqrt{\frac{p}{n}} g_\tau(\varepsilon_i) \eX^t_i \eu +D_i -\eE[D_i]\big].
\] 
Using assumption (A1) we have,  $ \eE\big[ (pn^{-1})^{1/2} g_\tau(\varepsilon_i) \eX^t_i \eu \big]=0$ and $ \Var\big[ ( pn^{-1})^{1/2} g_\tau(\varepsilon_i) \eX^t_i \eu \big]=p n^{-1} \eu^t  \sum^n_{i=1} \eX_i \eX_i^t \eu \Var [g_\tau(\varepsilon)] =O(p)$. Then, we have:
\begin{equation}
\label{D1}
B  \sqrt{\frac{p}{n}} g_\tau(\varepsilon_i) \eX^t_i \eu =O_{\PP}\big(p^{1/2} B \big).
\end{equation}
By Taylor expansion of $\rho_\tau \big(\varepsilon_i - B (pn^{-1})^{1/2} \eX_i^t \eu \big)$ around of $\varepsilon_i$, we can write $D_i$ also in the form: $D_i=2^{-1}B^2 p n^{-1}  |\eX_i^t \eu|^2  h_\tau(\widetilde \varepsilon_i)$, with $\widetilde \varepsilon_i$ a random variable between $\varepsilon_i$ and $\varepsilon_i+B (pn^{-1})^{1/2} \eX_i^t \eu$. On the other hand, given the definition of the function $h$, we have that:
\begin{equation}
\label{RR}
\PP \big[2 \min(\tau, 1-\tau) \leq h_\tau(\widetilde \varepsilon_i) \leq 2 \max(\tau, 1-\tau)\big]=1,
\end{equation}
and also $\Var [D_i] \leq \eE[D_i^2]=4^{-1}B^4 |\eX^t_i \eu|^4 \eE\big[h^2_\tau(\widetilde \varepsilon_i) \big]$. But $h^2_\tau(\widetilde \varepsilon_i)= 4\tau^2\e1_{\widetilde \varepsilon_i>0}+4(1-\tau)^2\e1_{\widetilde \varepsilon_i<0}$, then, $\eE\big[h^2_\tau(\widetilde \varepsilon_i) \big] \leq 4 \max \big(\tau^2, (1-\tau)^2 \big) \leq 4$. Thus, $\Var [D_i] \leq B^4   p^2 n^{-2} |\eX^t_i \eu|^4 $. On the other hand, the random variables $D_i$ defined by (\ref{Ri}), are independent. Then, 
\begin{equation}
\label{D2}
\sum^n_{i=1} \big[D_i- \eE[D_i]  \big]=O_{\PP}\bigg(\sqrt{\sum^n_{i=1} \Var [D_i]} \bigg) \leq O_{\PP} \bigg(\sqrt{\sum^n_{i=1} \eE[D_i^2]} \bigg) =O_{\PP}\bigg(B^2 \frac{p}{n^{1/2}}  \bigg).
\end{equation}
Relations (\ref{D1}) and (\ref{D2}), imply, since $p n^{-1} \rightarrow 0$, when $n \rightarrow \infty$, that:
\[
\Delta_1= O_{\PP}(B p^{1/2})+O_{\PP} (B^2 p n^{-1/2})= O_{\PP}(B p^{1/2}).
\]
Then, this last relation together relation (\ref{DD2}) imply  $\Delta_2 > | \Delta_1|$, with probability converging to one, for $B$ large enough. Relation (\ref{eq2}) follows, which implies the convergence rate of the expectile estimator.\\
\hh \textit{(ii)} For  $p$-vector $\eu=(u_1, \cdots , u_p)$, with $\| \eu\|_2=1$ and $B>0$ a constant, let us consider  the  difference 
\begin{equation} 
\label{RQ}
R_n\big( \ebo+B \sqrt{\frac{p}{n}} \eu\big)- R_n(\ebo)=Q_n\big( \ebo+B \sqrt{\frac{p}{n}} \eu\big)- Q_n(\ebo)+n \lambda_n \sum^p_{j=1} \widehat{\omega}_{n,j} \big[\big|\eb^0_j+B\sqrt{\frac{p}{n}}  u_j\big|- |\eb^0_j|  \big].
\end{equation}
The first term of the right-hand side of (\ref{RQ})  becomes by the above  proof  for \textit{(i)},
\begin{equation}
\label{QQ}
Q_n\big( \ebo+B \sqrt{\frac{p}{n}} \eu\big)- Q_n(\ebo)=B^2 \frac{p}{n} \eE \big[ h_\tau(\varepsilon) \big]  \sum^n_{i=1} (\eX_i^t \eu)^2 \big(1+o_{\PP}(1)\big)=O_{\PP}(B^2 p).
\end{equation}
Furthermore, for the penalty of (\ref{RQ}) we have:
\[
n \lambda_n \sum^p_{j=1} \widehat{\omega}_{n,j} \big[\big|\eb^0_j+B\sqrt{\frac{p}{n}}  u_j\big|- |\eb^0_j|  \big] \geq n \lambda_n \sum^{p_0}_{j=1} \widehat{\omega}_{n,j} \big[\big|\eb^0_j+B\sqrt{\frac{p}{n}}  u_j\big|- |\eb^0_j|  \big] \geq  - n \lambda_n \sum^{p_0}_{j=1} \widehat{\omega}_{n,j} B \sqrt{\frac{p}{n}} |u_j| .
\]
By the Cauchy-Schwarz inequality and afterwards by \textit{(i)}, we have
\begin{equation}
\label{PP}
\geq - n \lambda_n B \sqrt{\frac{p}{n}} \bigg( \sum^{p_0}_{j=1} \widehat{\omega}_{n,j}^2 \bigg)^{1/2} \| \eu \|_2= - B  \sqrt{\frac{p}{n}} n p_0^{1/2} \lambda_n= - B p_0^{1/2} n^{(1+c)/2} \lambda_n.
\end{equation}
Since $ p_0^{1/2} n^{(1-c)/2} \lambda_n \rightarrow 0$, as $n \rightarrow \infty$, we obtain that relation (\ref{QQ}) dominates (\ref{PP}) and the assertion regarding the   convergence rate of $\widehat{\eb}_n$ results.
\hspace*{\fill}$\blacksquare$  \\

\noindent {\bf Proof of Theorem \ref{Theorem 2SPL}}.\\
\textit{(i)}  Let us consider the parameter set:  ${\cal V}_p(\ebo) \equiv \big\{ \eb \in \R^p; \| \eb-\ebo\|_2 \leq B \sqrt{\frac{p}{n}}\big\}$, with $B>0$ large enough and ${\cal W}_n \equiv \acc{\eb \in {\cal V}_p(\ebo) ; \|\eb_{{\cal A}^c}\|_2 >0}$. According to Theorem \ref{th_vconv}, the estimator $\widehat{\eb}_n$ belongs to the set ${\cal V}_p(\ebo)$ with a probability converging to 1 as $n \rightarrow \infty$. In order to show the sparsity property of claim \textit{(i)}, we will show that, $\lim_{n \rightarrow \infty} \PP\big[\widehat{\eb}_n \in  {\cal W}_n \big]=0$. Note that if $\eb \in {\cal W}_n$, then $p >p_0$. \\ 
Let us consider two parameter vectors: $\eb=(\eb_{\cal A}, \eb_{{\cal A}^c}) \in {\cal W}_n$ and $ {\eb}^{(1)}=({\eb}_{\cal A}^{(1)}, {\eb}_{{\cal A}^c}^{(1)}) \in {\cal V}_p(\ebo)$,   such that  ${\eb}^{(1)}_{\cal A}=\eb_{\cal A}$ and $ {\eb}_{{\cal A}^c}^{(1)}=\textbf{0}_{p-p_0}$. For this parameters, we will study the following difference:
\begin{equation}
\label{DD}
 n^{-1} \big[R_n(\eb) - R_n({\eb}^{(1)}) \big] = n^{-1}\big[Q_n(\eb) - Q_n({\eb}^{(1)}) \big] +\lambda_n  \sum^p_{j=p_0+1} \widehat{\omega}_{n,j} | \beta_j|.
\end{equation}
For the first term of the right-hand side of   relation  (\ref{DD}), we have
\begin{align*}
&n^{-1} \sum^n_{i=1} \big[\rho_\tau(Y_i-\eX_i^t \eb^{(1)})- \rho_\tau(Y_i-\eX_i^t \eb)\big]
\\
 & \quad =  n^{-1} \sum^n_{i=1} \big[\rho_\tau\big(\varepsilon_i- \eX_{i, {\cal A}}^t (\eb_{\cal A}- \eb^0_{\cal A}) \big) -\rho_\tau \big( \varepsilon_i- \eX_{i, {\cal A}}^t (\eb_{\cal A}- \eb^0_{\cal A})-\eX^t_{i, {\cal A}^c} \eb_{{\cal A}^c} \big) \big]\\
  &\quad = n^{-1} \sum^n_{i=1} \big[ g(\varepsilon_i) \eX^t_{i,{\cal A} }\big(\eb_{\cal A}- \eb^0_{\cal A} \big)+ \frac{h(\varepsilon_i)}{2}\big(\eX^t_{i,{\cal A} }(\eb_{\cal A}- \eb^0_{\cal A} )\big)^2 +o_{\PP}\big(\eX^t_{i,{\cal A} }(\eb_{\cal A}- \eb^0_{\cal A} )\big)^2 \big] \\
  & \quad   -n^{-1} \sum^n_{i=1} \big[ g(\varepsilon_i) \eX^t_i\big(\eb_{\cal A}- \eb^0_{\cal A},\eb_{{\cal A}^c}\big)+ \frac{h(\varepsilon_i)}{2} \big(\eX^t_i (\eb_{\cal A}- \eb^0_{\cal A}, \eb_{{\cal A}^c}) \big)^2+o_{\PP}\big(\eX^t_i (\eb_{\cal A}- \eb^0_{\cal A}, \eb_{{\cal A}^c}) \big)^2\big] .
\end{align*}
By similar arguments used in  the proof of Theorem \ref{th_vconv}\textit{(i)} we have
\begin{align*}
&n^{-1} \sum^n_{i=1} g(\varepsilon_i) \eX^t_{i,{\cal A} }\big(\eb_{\cal A}- \eb^0_{\cal A} \big)  =\eE\big[ n^{-1} \sum^n_{i=1} g(\varepsilon_i) \eX^t_{i,{\cal A} }\big(\eb_{\cal A}- \eb^0_{\cal A} \big)\big] \\
 &  \hspace{6cm} +O_{\PP} \bigg( \Var \big[ n^{-1} \sum^n_{i=1} g(\varepsilon_i) \eX^t_{i,{\cal A} }\big(\eb_{\cal A}- \eb^0_{\cal A} \big)\big]\bigg)^{1/2}\\
&=O_{\PP} \bigg( \Var \big[ n^{-1} \sum^n_{i=1} g(\varepsilon_i) \eX^t_{i,{\cal A} }\big(\eb_{\cal A}- \eb^0_{\cal A} \big)\big]\bigg)^{1/2}   \leq O_{\PP} \bigg(\eE[g^2(\varepsilon_i)] \frac{1}{n^2} \sum^n_{i=1} \|\eX_{i,{\cal A} }\|^2_2 \|\eb_{\cal A}- \eb^0_{\cal A}\|^2_2 \bigg)^{1/2}  \\
&    =O_{\PP} \bigg( \frac{1}{n} \frac{p}{n}\bigg)^{1/2}=O_{\PP}\bigg(\frac{p^{1/2}}{n} \bigg).
\end{align*}
Proceeding similarly  as above, we get:
\[
n^{-1} \sum^n_{i=1} g(\varepsilon_i)\eX^t_i\big(\eb_{\cal A}- \eb^0_{\cal A},\eb_{{\cal A}^c}\big)=O_{\PP}\bigg(\frac{p^{1/2}}{n} \bigg).
\]
Taking into account  relation (\ref{RR}), we deduce that:
$
0< n^{-1} \sum^n_{i=1}  h(\varepsilon_i)  \big(\eX^t_i (\eb_{\cal A}- \eb^0_{\cal A}, \eb_{{\cal A}^c}) \big)^2=O_{\PP}(\| \eb_{\cal A}- \eb^0_{\cal A}\|^2_2 )=O_{\PP}( p n^{-1} )$ 
and also that,
\[
0< n^{-1} \sum^n_{i=1}  h(\varepsilon_{i}) \big(\eX^t_{i,{\cal A}} (\eb_{\cal A}- \eb^0_{\cal A}) \big)^2=O_{\PP}\bigg(\frac{p}{n} \bigg).
\]
By these relations, we obtain that the first term of the right-hand side of relation (\ref{DD}) is of order $p n^{-1}$. 
For the penalty of the right-hand side of relation (\ref{DD}), taking into account Theorem \ref{th_vconv}\textit{(i)} and since $\eb \in {\cal W}_n$ we obtain:
\[
\lambda_n  \sum^p_{j=p_0+1} \widehat{\omega}_{n,j} | \beta_j| \geq C\lambda_n \bigg( \frac{p}{n}\bigg)^{(1-\gamma)/2}.
\]
Using the supposition $\lambda_n (pn^{-1})^{-(1+\gamma)/2} {\underset{n \rightarrow \infty}{\longrightarrow}} \infty$, that is $\lambda_n n^{(1-c)(1+\gamma)/2} {\underset{n \rightarrow \infty}{\longrightarrow}} \infty$, we have that in the right-hand side of relation  (\ref{DD}), it's the penalty that dominates. Then, since  $n^{-1}\big[Q_n(\eb) - Q_n({\eb}^{(1)}) \big] =O_\PP(p n^{-1})$, we have,
\begin{equation}
\label{Rb}
n^{-1} \big[R_n(\eb) - R_n({\eb}^{(1)}) \big] \geq C\lambda_n \bigg( \frac{p}{n}\bigg)^{(1-\gamma)/2} .
\end{equation}
But, on the other hand, since $\eb^{(1)}_{{\cal A}^c}=\textbf{0}_{p-p_0}$, by similar arguments as above, we have,   $n^{-1} \big[R_n(\ebo) - R_n({\eb}^{(1)}) \big]=O_{\PP}(pn^{-1})$. From the last  relation together relation (\ref{Rb}), since $\lambda_n (p n^{-1})^{-(1+\gamma)/2} {\underset{n \rightarrow \infty}{\longrightarrow}} \infty$, we deduce, $\lim_{n  \rightarrow \infty}\PP [\widehat \eb_n \in {\cal W}_n] = 0$.\\
\hh \textit{(ii)} Given the previous result we consider the parameter vector  $\eb$ of the form: $\eb=\ebo+(p n^{-1})^{1/2} \ed$, with $\ed=(\ed_{\cal A}, \ed_{{\cal A}^c})$, $\ed_{{\cal A}^c}=\textbf{0}_{p-p_0}$, $\|\ed_{\cal A}\|_2 \leq C $. We study then the following difference:
\begin{equation}
\label{eq40}
\frac{1}{n}R_n\big(\ebo+\sqrt{\frac{p}{n}} \ed  \big)- \frac{1}{n}R_n\big(\ebo \big)=\frac{1}{n} \sum^n_{i=1}  \big[\rho_\tau\big(Y_i-\eX^t_i (\ebo+\sqrt{\frac{p}{n}} \ed ) \big)-\rho_\tau(\varepsilon_i) \big]+{\cal P}.
\end{equation}
For the  penalty ${\cal P}=\lambda_n \sum^{p_0}_{j=1} \widehat{\omega}_{n,j}\big(|\beta_j| -|\beta^0_j|\big) $ of the right-hand side of  relation (\ref{eq40}) we have, by \ref{th_vconv}\textit{(i)},  $\widehat{\omega}_{n,j}=|\widetilde{\beta}_{n,j}|^{-\gamma}=O_{\PP}(1)$ and  by the triangular inequality $\big||\beta_j| -|\beta^0_j|  \big| \leq |\beta_j -\beta^0_j|$. Then, as in the proof of Theorem \ref{th_vconv}, by relation (\ref{PP}), we obtain: 
\begin{equation}
\label{gp}
{\cal P}=O_{\PP}\bigg( \lambda_n  p_0^{1/2}  \bigg( \frac{p}{n}\bigg)^{1/2}\bigg)=O_{\PP}\bigg( \lambda_n p_0^{1/2} n^{(c-1)/2 } \bigg).
\end{equation}
For the  first term of the right-hand side of relation (\ref{eq40}) we have:
\begin{align*}
&\frac{1}{n} \sum^n_{i=1}  \big[\rho_\tau\big(Y_i-\eX^t_{i,{\cal A}} (\eb^0_{\cal A}+\sqrt{\frac{p}{n}} \ed _{\cal A}) \big)-\rho_\tau(\varepsilon_i) \big] \nonumber
\\
& \qquad =-\frac{1}{n} \sum^n_{i=1}  g(\varepsilon_i) \big(\eX^t_{i,{\cal A}}\ed _{\cal A} \big)\sqrt{\frac{p}{n}}+\frac{1}{2n}  \sum^n_{i=1} \big[\frac{p}{n}\|\eX^t_{i,{\cal A}}\ed _{\cal A}\|^2_2 h(\varepsilon_i)  +o_{\PP} \big(  \|\eX^t_{i,{\cal A}}\ed _{\cal A}\|^2\big) \big] 
\end{align*}
\begin{align}
& \qquad=\bigg(-\frac{1}{n}\sqrt{\frac{p}{n}} \sum^n_{i=1}  g(\varepsilon_i) \big(\eX^t_{i,{\cal A}}\ed _{\cal A} \big)+\frac{1}{2n} \frac{p}{n} \sum^n_{i=1}  \big(\ed _{\cal A}^t \eX_{i,{\cal A}}\eX^t_{i,{\cal A}}\ed _{\cal A} h(\varepsilon_i) \big)\bigg) \big( 1+o_{\PP}(1)\big) \nonumber\\
&\qquad =\bigg(-\frac{1}{n}\sqrt{\frac{p}{n}} \sum^n_{i=1}  g(\varepsilon_i) \big(\eX^t_{i,{\cal A}}\ed _{\cal A} \big) \nonumber \\
&\qquad  \qquad \qquad +\frac{1}{2n} \frac{p}{n} \sum^n_{i=1}  \big(\ed _{\cal A}^t \eX_{i,{\cal A}}\eX^t_{i,{\cal A}}\ed _{\cal A} \big(\eE[h(\varepsilon_i)]+h(\varepsilon_i)-\eE[h(\varepsilon_i)] \big)\big)\bigg) \big( 1+o_{\PP}(1)\big)\nonumber\\
& \qquad=\bigg(-\frac{1}{n}\sqrt{\frac{p}{n}} \sum^n_{i=1}  g(\varepsilon_i) \big(\eX^t_{i,{\cal A}}\ed _{\cal A} \big)+\frac{1}{2n} \frac{p}{n} \sum^n_{i=1}  \big(\ed _{\cal A}^t \eX_{i,{\cal A}}\eX^t_{i,{\cal A}}\ed _{\cal A}  \eE[h(\varepsilon_i)]\big)\bigg) \big( 1+o_{\PP}(1)\big),
\label{gg}
\end{align}
which has as minimizer the solution of
\[
-\frac{1}{n}\sqrt{\frac{p}{n}} \sum^n_{i=1}  g(\varepsilon_i) \eX_{i,{\cal A}}+\eU_{n,{\cal A}}\sqrt{\frac{p}{n}}  \ed _{\cal A} \eE[h(\varepsilon)]=\textbf{0}_{p_0},
\]
from where, we get, 
\[
\sqrt{\frac{p}{n}}  \ed _{\cal A} =\frac{\eU^{-1}_{n,{\cal A}}}{\eE[h(\varepsilon)]}\frac{1}{n} \sum^n_{i=1} g(\varepsilon_i) \eX_{i,{\cal A}}.
\]
We deduce that, the minimum value of (\ref{gg}) is of order $O_{\PP}\big(p n^{-1}\|\ed _{\cal A}\|_2 \big)=O_{\PP}\big(p n^{-1}\big)=O_{\PP}\big(n^{c-1}\big)$. Taking into account the supposition $\lambda_n p_0^{1/2} n^{(1-c)/2} {\underset{n \rightarrow \infty}{\longrightarrow}} 0$ and  relation (\ref{gp}), we have that, ${\cal P}=o_{\PP}(p n^{-1}) $. Then, in the right-hand side of relation (\ref{eq40}), the first term is the dominant one. \\
Let us now consider the following random variable sequence:
\[
W_i \equiv g(\varepsilon_i) \eu^t \frac{\eU^{-1}_{n,{\cal A}}}{\eE[h(\varepsilon)]}  \eX_{i,{\cal A}},
\]
with $\eu \in \R^{p_0}$, $\| \eu \|_2=1$. For the random variable   $W_i$, we have that $\eE[W_i]=0$ and $\Var[W_i]=\eE^{-2}[h(\varepsilon)]\eU^{-1}_{n,{\cal A}} \eu^t \eX_{i,{\cal A}} \eX^t_{i,{\cal A}} \eu \Var[g(\varepsilon_i) ]$. Thus, taking into account assumption (A1), we get: 
\[
\sum^n_{i=1}\Var[W_i]=n \frac{\eu^t \eU^{-1}_{n,{\cal A}} \eu}{\eE^2[h(\varepsilon)]}\Var[g(\varepsilon) ],
\]
which implies
\[
\sqrt{n} \frac{\eE[h(\varepsilon)]}{\sqrt{\Var[g(\varepsilon) ]}}\frac{\eu^t \big( \widehat \eb_n -\ebo\big)_{\cal A}}{\big(\eu^t \eU^{-1}_{n,{\cal A}} \eu \big)^{1/2} } \overset{\cal L} {\underset{n \rightarrow \infty}{\longrightarrow}} {\cal N}(0,1).
\]
The proof of claim \textit{(ii)} is finished.
\hspace*{\fill}$\blacksquare$  \\

  \subsection{\textbf{Result proofs  in Section \ref{casec>1}}}
  \noindent {\bf Proof of Theorem \ref{th_vconvn}}.\\
The proof is similar to that of Theorem \ref{th_vconv}. Consequently,, we give only the main results.  Otherwise, instead of the Cauchy-Schwarz inequality we use Holder's inequality: $|\eX_i^t \eu| \leq \|\eX_i\|_{\infty} \|\eu\|_1$ and then we obtain: $0 <\Delta_2=O(B^2 b^2_n n \|\eu\|^2_1)$.  \\
For a  $p$-vector $\eu=(u_1, \cdots , u_p)$, with $\| \eu\|_1=1$ and  a constant $B>0$, let be the difference
\begin{equation} 
\label{RQn}
R_n\big( \ebo+B b_n \eu\big)- R_n(\ebo)=Q_n\big( \ebo+Bb_n \eu\big)- Q_n(\ebo)+n \lambda_n \sum^p_{j=1} \widehat{\omega}_{n,j} \big[\big|\eb^0_j+B b_n  u_j\big|- |\eb^0_j|  \big].
\end{equation}
By a similar approach made for the terms $\Delta_1$ and $\Delta_2$ of the proof of Theorem \ref{th_vconv}, we  obtain:
\begin{equation}
\label{QQn}
Q_n\big( \ebo+B b_n \eu\big)- Q_n(\ebo)=O_{\PP}(B^2 b_n^2 n\|\eu\|^2_1).
\end{equation}
For the penalty of the right-hand side of relation (\ref{RQn}) we have:
\[
n \lambda_n \sum^p_{j=1} \widehat{\omega}_{n,j} \big[\big|\eb^0_j+Bb_n  u_j\big|- |\eb^0_j|  \big] \geq n \lambda_n \sum^{p_0}_{j=1} \widehat{\omega}_{n,j} \big[\big|\eb^0_j+Bb_n  u_j\big|- |\eb^0_j|  \big] \geq  - n \lambda_n \sum^{p_0}_{j=1} \widehat{\omega}_{n,j} B b_n |u_j| ,
\]
by the Cauchy-Schwarz inequality and afterwards by the estimator consistency of $\overset{\vee}{{\eb}}_n$, we have 
\begin{equation}
\label{PPn}
\geq - n \lambda_n B b_n \bigg( \sum^{p_0}_{j=1} \widehat{\omega}_{n,j}^2 \bigg)^{1/2} \| \eu \|_2\geq  - B C b_n n p_0^{1/2} \lambda_n\|\eu\|^2_1= - B n \lambda_np_0^{1/2} b_n.
\end{equation}
Since $ \lambda_n p_0^{1/2} b^{-1}_n \rightarrow 0$, as $n \rightarrow \infty$, then relation (\ref{QQn}) dominates (\ref{PPn}) and the theorem follows.
\hspace*{\fill}$\blacksquare$  \\

\noindent {\bf Proof of Theorem \ref{Theorem 2SPLn}}.\\
\textit{(i)} Let $j \in {\cal A}^c$ be, then  $j >p_0$. Thus, the derivative of the random process $R_n(\eb)$ in respect to  $\beta_j$ is:
\begin{equation}
\frac{\partial R_n(\eb)}{\partial \beta_j} = \sum^n_{i=1} g_\tau\big(Y_i- \eX^t_i\eb \big) X_{ij}+n \lambda_n \widehat{\omega}_{n,j} \textrm{sgn}(\beta_j).
 \label{tt}
\end{equation}
For the first term of the right-hand side of relation (\ref{tt}), we have, $ g_\tau\big(Y_i- \eX^t_i\eb \big)=g_\tau\big(\varepsilon_i -\eX_i^t(\eb-\ebo)  \big)$. We denote $\eta_i=\eX_i^t(\eb-\ebo) $ and then $g_\tau(\varepsilon_i-\eta_i)=g_\tau(\varepsilon_i)-\eta_ih_\tau(\tilde \eta_i)$, with $\widetilde{\eta}_i$ a random variable variable between  $\varepsilon_i$ and $\varepsilon_i-\eta_i$.  Then
\[
\sum^n_{i=1} g_\tau\big(Y_i- \eX^t_i\eb \big) X_{ij}=\sum^n_{i=1} g_\tau(\varepsilon_i) X_{ij}- \sum^n_{i=1} \eX_i^t(\eb-\ebo)  h_\tau(\widetilde \eta_i) X_{ij}.
\]
By the Central Limit Theorem, tacking into account assumption  (A4), we have that:  $\sum^n_{i=1} g_\tau(\varepsilon_i) X_{ij}=O_{\PP}(n^{1/2})$. On the other hand, $0 < h_\tau(\widetilde \eta_i) <2$ with probability 1. Using the   Holder's inequality, we have, with probability one, 
\[
\bigg|  \sum^n_{i=1} \eX_i^t(\eb-\ebo)  h_\tau(\widetilde \eta_i) X_{ij} \bigg| \leq \sum^n_{i=1}  \bigg|  \eX_i^t(\eb-\ebo)  h_\tau(\widetilde \eta_i) X_{ij}\bigg| \leq \sum^n_{i=1}\| \eX_i\|_{\infty} \big\|h_\tau(\widetilde \eta_i) X_{ij}(\eb-\ebo)\big\|_1 ,
\]
from where, tacking into account assumption (A4), we have:  $\sum^n_{i=1} \eX_i^t(\eb-\ebo)  h_\tau(\widetilde \eta_i) X_{ij}=O_{\PP}(nb_n)$. Thus,
\begin{equation}
\label{t1}
\sum^n_{i=1} g_\tau\big(Y_i- \eX^t_i\eb \big) X_{ij}=O_{\PP}(n b_n).
\end{equation}
For the penalty of relation (\ref{tt}) we have: $n \lambda_n \widehat{\omega}_{n,j} =O_{\PP} \big(n \lambda_n \min \big(n^{1/2} , a_n^{-\gamma}\big)\big)$. Since, $\lambda_n b_n^{-1}\min \big(n^{1/2}, a_n^{-\gamma}\big) \rightarrow \infty$, as $n \rightarrow \infty$, tacking also into account  relation (\ref{t1}) we have that:
\[
\frac{\partial R_n(\eb)}{\partial \beta_j} \left\{
\begin{array}{lll}
>0,& &\textrm{if } \beta_j>0,\\
& & \\
<0, & & \textrm{if } \beta_j<0.
\end{array}
\right.
\]
The function $R_n(\eb)$ is continuous in $\eb$. Then, the solution of (\ref{tt}) must be equal to 0.  From where  $\widehat{\eb}_{n,{\cal A}^c}=\textbf{0}_{p-p_0}$, with probability converging to 1. This relation implies $\widehat{{\cal A}}_n \subseteq {\cal A}$ with probability converging to 1 when $n \rightarrow \infty$.\\
 On the basis of this result, from now on we consider the parameters $\eb$ of the form $\eb=\big(\eb_{\cal A}, \textbf{0}_{p-p_0}  \big)$. We must show now that ${\cal A} \subseteq \widehat{{\cal A}}_n $. By Theorem \ref{th_vconvn} we have $\| \widehat{\eb}_{\cal A} -\eb^0_{\cal A} \|_1=O_{\PP}\big(b_n \big)$, from where for any  $j=1, \cdots , p_0$, we obtain, $\widehat{\beta}_{n,j} \overset{\PP} {\underset{n \rightarrow \infty}{\longrightarrow}}  \beta^0_j \neq 0$. Thus, since $b_n {\underset{n \rightarrow \infty}{\longrightarrow}}  0$, we have that $\widehat{\beta}_{n,j} \neq 0$ with probability converging to 1, from where ${\cal A} \subseteq \widehat{{\cal A}}_n $.\\
\hh \textit{(ii)} Given the previous result \textit{(i)} and Theorem \ref{th_vconvn}, we consider the parameters  $\eb$ of the form: $\eb=\ebo+b_n \ed$, with $\ed=(\ed_{\cal A}, \ed_{{\cal A}^c})$, $\ed_{{\cal A}^c}=\textbf{0}_{p-p_0}$, $\|\ed_{\cal A}\|_1 \leq C $.
For the penalty ${\cal P}$ of the right-hand side of   relation (\ref{eq40}) we have:  $\big|{\cal P} \big|=\lambda_n \big|\sum^{p_0}_{j=1} \widehat{\omega}_{n,j}|\beta_j| -|\beta^0_j|] \big|  \leq \lambda_n \sum^{p_0}_{j=1} \widehat{\omega}_{n,j} \big| \beta_j -\beta^0_j \big| \leq \lambda_n \big( \sum^{p_0}_{j=1} \widehat{\omega}_{n,j}^2\big)^{1/2} \| (\eb-\ebo)_{\cal A}\|_2=O_{\PP} \big(\lambda_n b_n p^{1/2}_0 \big)$.  
For the main part, we have:
\[
\frac{1}{n} \sum^n_{i=1}  \big[\rho_\tau\big(Y_i-\eX^t_{i,{\cal A}} (\eb^0_{\cal A}+b_n \ed _{\cal A}) \big)-\rho_\tau(\varepsilon_i) \big] 
\]
\[
=\bigg(-\frac{1}{n}b_n \sum^n_{i=1}  g(\varepsilon_i) \big(\eX^t_{i,{\cal A}}\ed _{\cal A} \big)+\frac{1}{2n} b^2_n \sum^n_{i=1}  \big(\ed _{\cal A}^t \eX_{i,{\cal A}}\eX^t_{i,{\cal A}}\ed _{\cal A}  \eE[h(\varepsilon_i)]\big)\bigg) \big( 1+o_{\PP}(1)\big).
\]
The end of the proof is similar to that of Theorem \ref{th_vconvn}\textit{(ii)}. 
\hspace*{\fill}$\blacksquare$  \\


\end{document}